# Pascal Pyramids, Pascal Hyper-Pyramids and a Bilateral Multinomial Theorem

Martin Erik Horn, University of Potsdam
Am Neuen Palais 10, D - 14469 Potsdam, Germany
E-Mail: *marhorn@rz.uni-potsdam.de*

**Abstract**
Part I: The two-dimensional Pascal Triangle will be generalized into a three-dimensional Pascal Pyramid and four-, five- or whatsoever-dimensional hyper-pyramids.
Part II: The Bilateral Binomial Theorem will be generalised into a Bilateral Trinomial Theorem resp. a Bilateral Multinomial Theorem.

**Introduction**

The complete Pascal Plane with its three Pascal Triangles consists of the following numbers

$$f_{(x, y)} = \lim_{h \to 0} \frac{(x+y+h)!}{(x+h)! \cdot (y+h)!} = \lim_{h \to 0} \frac{\Gamma(x+y+1+h)}{\Gamma(x+1+h) \cdot \Gamma(y+1+h)} \quad (1)$$

and looks like this if the positive directions are pointed downwards:

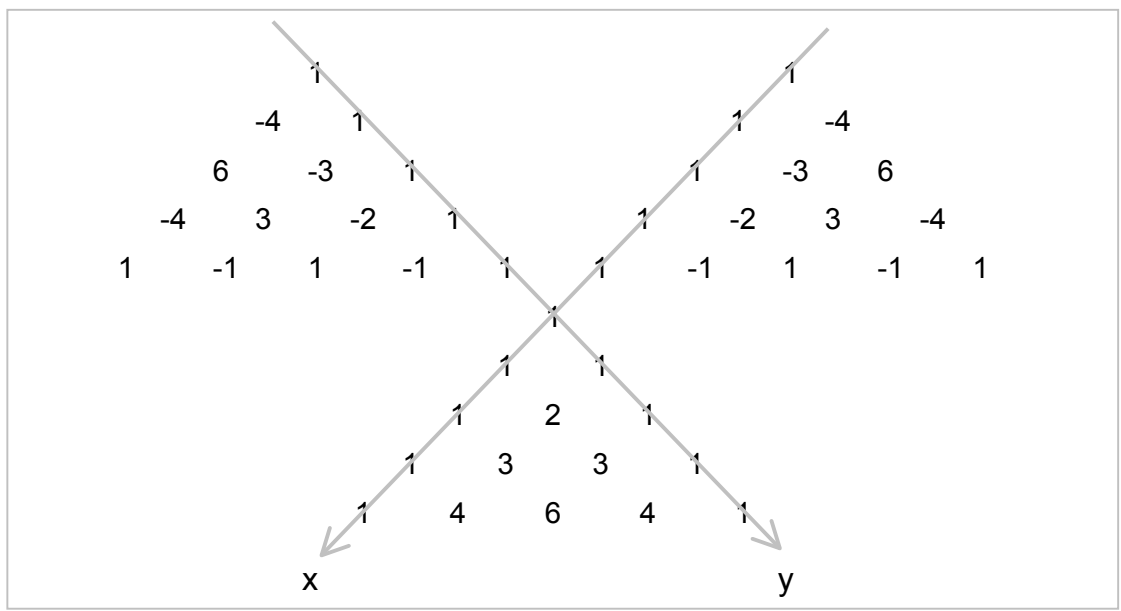

If now numbers with distance 1 are added and the definition of the bilateral hypergeometric function of [1] is used

$$_1H_1(a; b; z) = ... + \frac{(b-1) \cdot (b-2)}{(a-1) \cdot (a-2)} \cdot z^{-2} + \frac{b-1}{a-1} \cdot z^{-1} + 1 + \frac{a}{b} \cdot z + \frac{a \cdot (a+1)}{b \cdot (b+1)} \cdot z^2 + ... \quad (2)$$

the following bilateral hypergeometric identity is reached:

$$2^x = \frac{x!}{y! \cdot (x-y)!} \cdot {}_1H_1[(y-x); (y+1); -1] \qquad x, y \in \mathbf{R} \quad (3)$$

This is a special case of the Bilateral Binomial Theorem [2, 3] with $|z| = 1$:

$$(1+z)^x = \frac{\Gamma(x+1)}{\Gamma(y+1) \cdot \Gamma(x-y+1)} \cdot {}_1H_1[(y-x); (y+1); -z] \qquad x, y \in \mathbf{R}; \; z \in \mathbf{C} \quad (4)$$

**Part I: Pascal Pyramids and Pascal Hyper-Pyramids**

The Pascal Plane, which consists of binomial coefficients, can be generalized into the Pascal Space using trinomial coefficients

$$(x_1, x_2, x_3) = \frac{(x_1 + x_2 + x_3)!}{x_1! \cdot x_2! \cdot x_3!} \tag{5}$$

Then the Pascal Pyramid can be constructed by adding every three appropriate neighbouring numbers and writing the result beneath them:

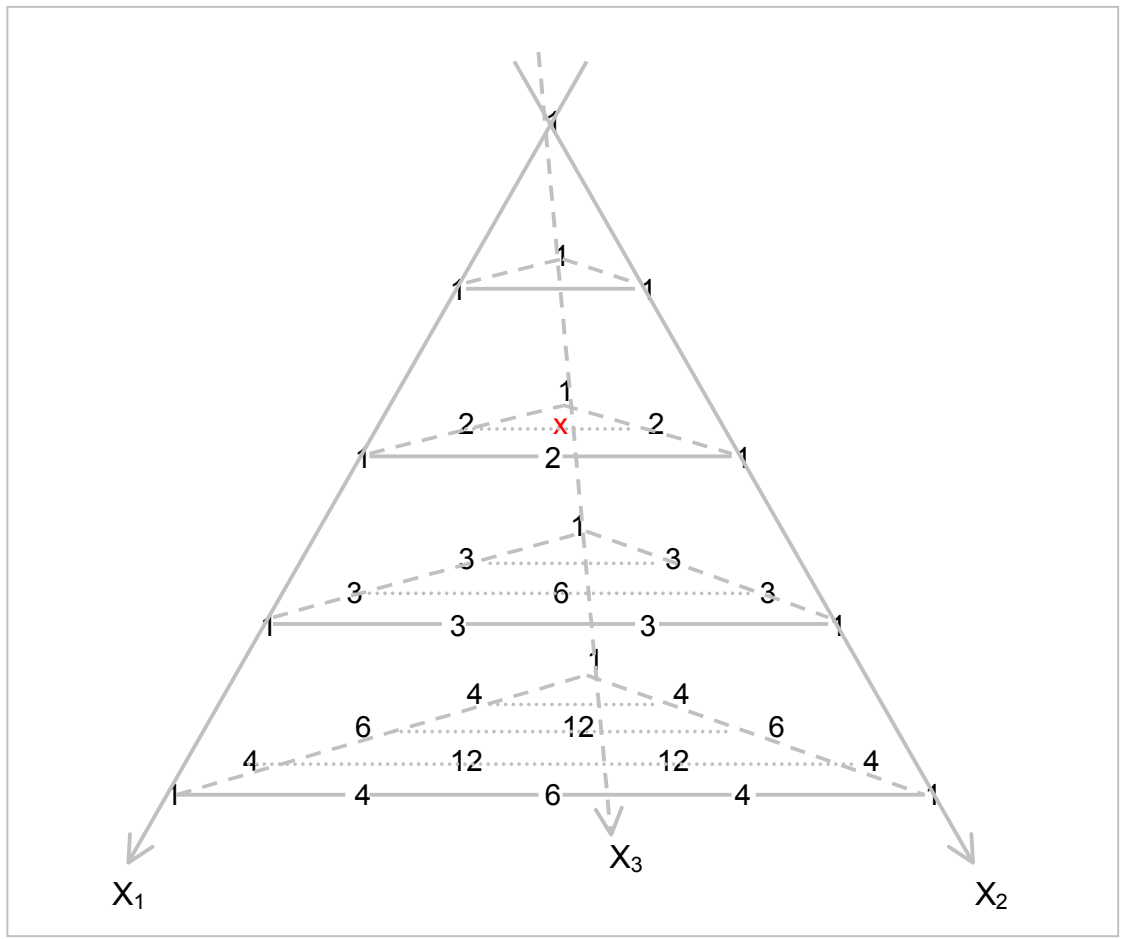

Remark: No, there isn't a proud 3 sitting in the middle of the second triangle at the marked red position x. This is the place for the following humble trinomial coefficient:

$$\left(\tfrac{2}{3}, \tfrac{2}{3}, \tfrac{2}{3}\right) = \frac{2!}{\tfrac{2}{3}! \cdot \tfrac{2}{3}! \cdot \tfrac{2}{3}!} = 3 \cdot \left(-\tfrac{1}{3}, \tfrac{2}{3}, \tfrac{2}{3}\right) \tag{6}$$

because the construction law of trinomial coefficients reads:

$$(x_1, x_2, x_3) = (x_1 - 1, x_2, x_3) + (x_1, x_2 - 1, x_3) + (x_1, x_2, x_3 - 1) \tag{7}$$

But the picture above shows only a quarter of the truth, of course, for three similar pyramids can be constructed in the negative coordinate region using these numbers

$$f_{(x, y, z)} = \lim_{h \to 0} \frac{(x + y + z + h)!}{(x + h)! \cdot (y + h)! \cdot (z + h)!} = \lim_{h \to 0} \frac{\Gamma(x + y + z + 1 + h)}{\Gamma(x + 1 + h) \cdot \Gamma(y + 1 + h) \cdot \Gamma(z + 1 + h)} \tag{8}$$

as the following drawing indicates:

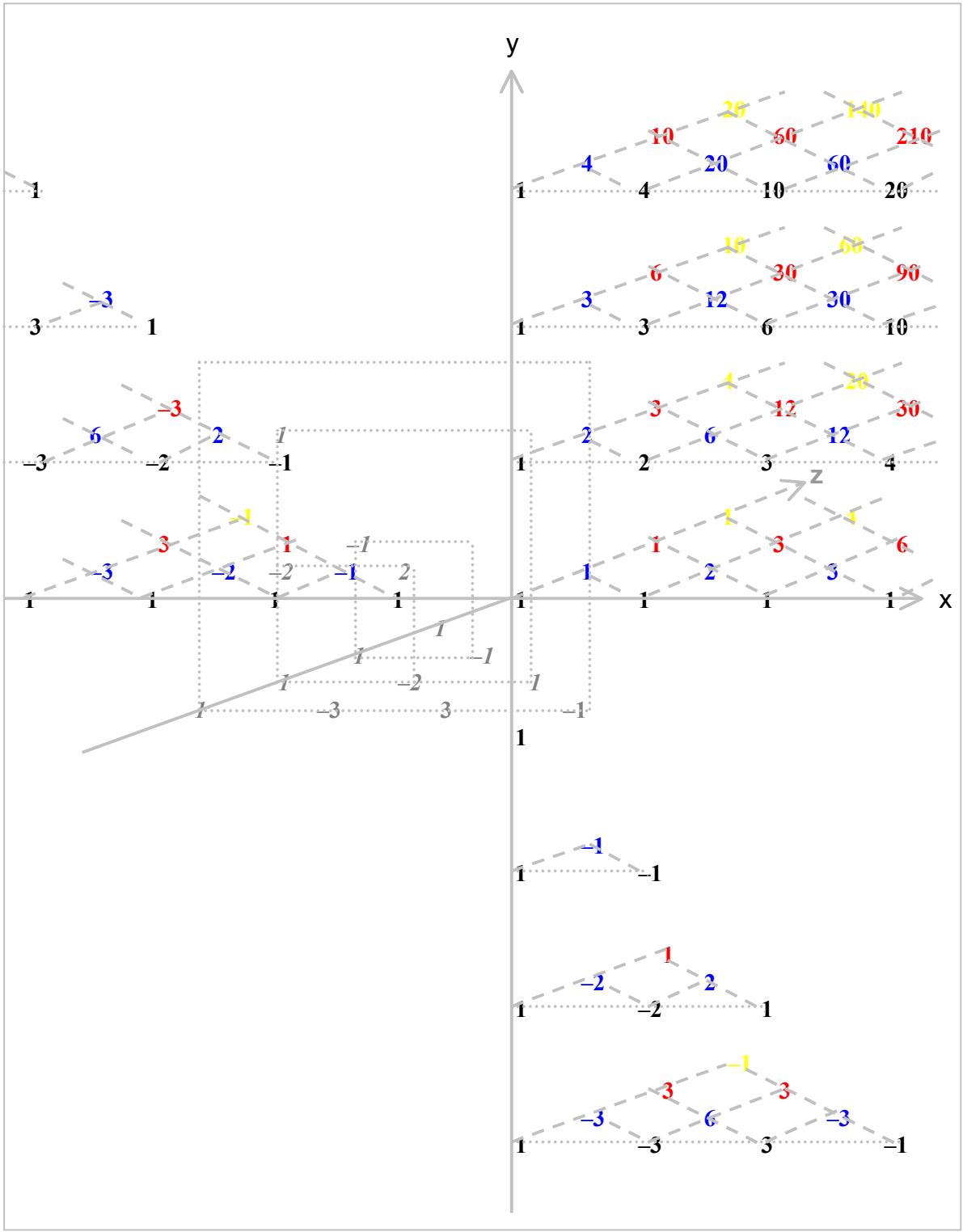

And slight rotations of the axes produce a more symmetric design with tetrahedral order as the picture on the right is supposed to show.

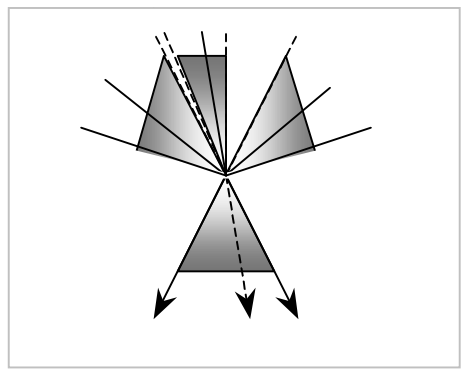

The next step is to increase the dimension again by considering quatronomial coefficients, which fill the four-dimensional Pascal Hyper-Space:

$$(x_1, x_2, x_3, x_4) = \frac{(x_1 + x_2 + x_3 + x_4)!}{x_1! \cdot x_2! \cdot x_3! \cdot x_4!} \tag{9}$$

$$\begin{aligned} = (x_1 - 1, x_2, x_3, x_4) + (x_1, x_2 - 1, x_3, x_4) \\ + (x_1, x_2, x_3 - 1, x_4) + (x_1, x_2, x_3, x_4 - 1) \end{aligned} \tag{10}$$

By again using

$$f_{(w, x, y, z)} = \lim_{h \to 0} \frac{\Gamma(w + x + y + z + 1 + h)}{\Gamma(w + 1 + h) \cdot \Gamma(x + 1 + h) \cdot \Gamma(y + 1 + h) \cdot \Gamma(z + 1 + h)} \tag{11}$$

five Pascal Hyper-Pyramids can be found. The three-dimensional hyper-surfaces of these four-dimensional hyper-pyramids consist of the Pascal Pyramids (with some more minus-signs every now and then) sketched on the previous page.

This procedure can be continued till eternity. The multinomial coefficients

$$(x_1, x_2, ..., x_n) = \frac{(x_1 + x_2 + ... + x_n)!}{x_1! \cdot x_2! \cdot ... \cdot x_n!} \tag{12}$$

$$= (x_1 - 1, x_2, ..., x_4) + (x_1, x_2 - 1, ..., x_4) + ... + (x_1, x_2, ..., x_4 - 1) \tag{13}$$

live in n-dimensional Pascal Hyper-Space, and with the help of

$$f_{(x_1, x_2, ..., x_n)} = \lim_{h \to 0} \frac{\Gamma(x_1 + x_2 + ... + x_n + 1 + h)}{\Gamma(x_1 + 1 + h) \cdot \Gamma(x_2 + 1 + h) \cdot ... \cdot \Gamma(x_n + 1 + h)} \tag{14}$$

n+1 Pascal Hyper-Pyramids can be constructed. These n-dimensional hyper-pyramids possess (n − 1)-dimensional hyper-surfaces which look like the Pascal Pyramids of one dimension less and some more minus-signs every now and then.

**Part II: Bilateral Multinomial Theorems**

Formula (3) was found by adding numbers of distance 1 which lie on a straight line in the Pascal Plane. One dimension higher a similar formula should be found, if all numbers of distance 1 of the Pascal Space are added which lie in a straight plane. This then would result in powers of 3

$$3^n = \sum_{y=-\infty}^{\infty} \sum_{x=-\infty}^{\infty} (x; y; n - x - y) \qquad x, y \in \mathbf{R} \tag{15}$$

if the series converged. But this double bilateral summation isn't supposed to converge for it is a special case ($|z_1| = |z_2| = 1$) of the Bilateral Trinomial Theorem

$$(1 + z_1 + z_2)^n = \sum_{x_2=-\infty}^{\infty} \sum_{x_1=-\infty}^{\infty} (x_1; x_2; n - x_1 - x_2) \cdot z_1^{x_1} \cdot z_2^{x_2} \qquad \begin{aligned} x_1, x_2 &\in \mathbf{R} \\ z_1, z_2 &\in \mathbf{C} \end{aligned} \tag{16}$$

The Bilateral Trinomial Theorem can be reformulated as

$$(1 + x + y)^n = \sum_{\ell=-\infty}^{\infty} \sum_{k=-\infty}^{\infty} \frac{x^k \cdot y^\ell}{k! \cdot \ell! \cdot (n+1)_{-k-\ell}} \tag{17}$$

where $(a)_k$ denotes the Pochhammer Symbol

$$(a)_k = \frac{\Gamma(a+k)}{\Gamma(a)} \qquad \text{resp.} \qquad (a+k)_{-k} = \frac{\Gamma(a)}{\Gamma(a+k)} \qquad (18)$$

Using the results of [2, 3] the double summation can be evaluated easily, giving a proof of formula (17) for special values of x and y.

$$\sum_{\ell=-\infty}^{\infty}\sum_{k=-\infty}^{\infty} \frac{x^k \cdot y^\ell}{k! \cdot \ell! \cdot (n+1)_{-k-\ell}} = \sum_{\ell=-\infty}^{\infty}\left(\binom{n}{\ell} \cdot \sum_{k=-\infty}^{\infty}\binom{n-\ell}{k} \cdot x^k \cdot y^\ell\right) \qquad (19)$$

Of course the binomial coefficients of (19) are generalized here as

$$\binom{n}{\ell} = \lim_{h\to\infty} \frac{(n+h)!}{(\ell+h)! \cdot (n-\ell+h)!} = \lim_{h\to\infty} \frac{\Gamma(n+h-1)}{\Gamma(\ell+h-1)\cdot\Gamma(n-\ell+h-1)} \qquad (20)$$

With $|x|=1$ this gives

$$\sum_{\ell=-\infty}^{\infty}\sum_{k=-\infty}^{\infty} \frac{x^k \cdot y^\ell}{k! \cdot \ell! \cdot (n+1)_{-k-\ell}} = \sum_{\ell=-\infty}^{\infty}\left(\binom{n}{\ell} \cdot (1+x)^{n-\ell} \cdot y^\ell\right) \qquad (21)$$

$$= (1+x)^n \cdot \sum_{\ell=-\infty}^{\infty}\left(\binom{n}{\ell}\left(\frac{y}{1+x}\right)^\ell\right) \qquad (22)$$

And with $|y|=|1+x|$ the expected result emerges:

$$\sum_{\ell=-\infty}^{\infty}\sum_{k=-\infty}^{\infty} \frac{x^k \cdot y^\ell}{k! \cdot \ell! \cdot (n+1)_{-k-\ell}} = (1+x)^n \cdot \left(1+\frac{y}{1+x}\right)^n \qquad (23)$$

$$= (1+x+y)^n \qquad (24)$$

The same strategy leads to a Bilateral Quatronomial Theorem:

$$(1+x+y+z)^n = \sum_{m=-\infty}^{\infty}\sum_{\ell=-\infty}^{\infty}\sum_{k=-\infty}^{\infty} \frac{x^k \cdot y^\ell \cdot z^m}{k! \cdot \ell! \cdot m! \cdot (n+1)_{-k-\ell-m}} \qquad (25)$$

with $|x|=1$, $|y|=|1+x|$ and $|z|=|1+x+y|$.

And this again can be extended till eternity giving the Bilateral Multinomial Theorem:

$$\boxed{(1+\sum_{i=1}^{k} x_i)^n = \sum_{\ell_1=-\infty}^{\infty}\sum_{\ell_2=-\infty}^{\infty}\cdots\sum_{\ell_k=-\infty}^{\infty} \frac{1}{(n+1)_{-\sum_{i=1}^{k}\ell_i}} \cdot \prod_{i=1}^{k} x_i^{\ell_i} \Big/ \prod_{i=1}^{k}(\ell_i!)} \qquad (26)$$

with $|x_i|=|1+\sum_{j=1}^{i-1} x_j|$ and $\begin{array}{l}\ell_i \in \mathbf{R} \\ x_i \in \mathbf{C}\end{array}$ .

**Epilogue**

To increase the aesthetical value of the indicated results a more symmetric formulation of the Bilateral Multinomial Theorem (26) can be given:

$$(x_0 + x_1 + x_2 + ... + x_k)^n$$

$$= \sum_{\ell_0=-\infty}^{\infty} \sum_{\ell_1=-\infty}^{\infty} \cdots \sum_{\ell_k=-\infty}^{\infty} \frac{\left(x_0 - \frac{1}{k+1}\right)^{\ell_0} \cdot \left(x_1 - \frac{1}{k+1}\right)^{\ell_1} \cdot ... \cdot \left(x_k - \frac{1}{k+1}\right)^{\ell_k}}{\ell_0! \cdot \ell_1! \cdot ... \cdot \ell_k! \cdot (n+1)_{-\ell_0-\ell_1-...-\ell_k}} \qquad (27)$$

But this of course doesn't change the fact that convergence is possible only with an unsymmetrical handling of the variables:

$$x_i = \left| 1 + x_1 - \frac{1}{k} + x_2 - \frac{1}{k} + ... + x_{i-1} - \frac{1}{k} \right| \qquad (28)$$